\documentclass[12pt,british,a4paper]{amsart}

\usepackage{lmodern}
\usepackage[T1]{fontenc}
\usepackage{microtype}

\usepackage{babel,amssymb,amscd,url,stmaryrd,mathtools,fixltx2e}
\usepackage{paralist,csquotes,stackrel}
\usepackage[mathscr]{eucal}

\numberwithin{equation}{section}
\mathtoolsset{mathic=true}
\setdefaultenum{\upshape (i)}{}{}{}

\theoremstyle{plain}
\newtheorem{theorem}[equation]{Theorem}

\theoremstyle{definition}

\theoremstyle{remark}
\newtheorem{Example}[equation]{Example}
\newtheorem{remark}[equation]{Remark}

\newenvironment{example}{\begin{Example}\pushQED{\qee}}{\popQED\end{Example}}
\DeclareRobustCommand{\qee}{%
  \ifmmode \mathqee
  \else
    \leavevmode\unskip\penalty9999 \hbox{}\nobreak\hfill
    \quad\hbox{\qeesymbol}%
  \fi
}
\newcommand{\mathqee}{\quad\hbox{\qeesymbol}}
\newcommand{\qeesymbol}{\ensuremath\diamondsuit}

\newcommand{\Lie}[1]{\operatorname{\textsl{#1}}}
\newcommand{\lie}[1]{\operatorname{\mathfrak{#1}}}
\newcommand{\GL}{\Lie{GL}}

\newcommand{\sln}{\lie{sl}}
\newcommand{\SO}{\Lie{SO}}
\newcommand{\so}{\lie{so}}
\newcommand{\SP}{\Lie{Sp}}
\newcommand{\sP}{\lie{sp}}

\newcommand{\bmf}{\lie b}

\newcommand{\kf}{\lie k}
\newcommand{\n}{\lie n}
\newcommand{\p}{\lie p}

\newcommand{\SL}{\Lie{SL}}
\newcommand{\SU}{\Lie{SU}}
\newcommand{\Un}{\Lie{U}}

\newcommand{\jj}{\mathbf j}

\newcommand{\C}{{\mathbb C}}
\newcommand{\HH}{{\mathbb H}}

\newcommand{\R}{{\mathbb R}}
\newcommand{\Z}{{\mathbb Z}}

\newcommand{\cO}{\mathcal O}
\newcommand{\tH}{\tilde H}

\newcommand{\hkimpl}{\textup{hkimpl}}
\newcommand{\impl}{\textup{impl}}

\newcommand{\symp}{{\sslash}} 
\newcommand{\hkq}{{\sslash\mkern-6mu/}}

\DeclareMathOperator{\Hom}{Hom}

\DeclareMathOperator{\im}{im}

\DeclareMathOperator{\Spec}{Spec}

\DeclarePairedDelimiter{\abs}{\lvert}{\rvert}

\begin{document}

\title{Symplectic and Hyperk\"ahler implosion}

\author[A.~Dancer]{Andrew Dancer}
\address[Dancer]{Jesus College\\
Oxford\\
OX1 3DW\\
United Kingdom} \email{dancer@maths.ox.ac.uk}

\author[B.~Doran]{Brent Doran}
\address[Doran]{Department of Mathematics\\
ETH Z\"urich\\
8092 Z\"urich\\
Switzerland}
\email{brent.doran@math.ethz.ch}

\author[F.~Kirwan]{Frances Kirwan}
\address[Kirwan]{Balliol College\\
Oxford\\
OX1 3BJ\\
United Kingdom} \email{kirwan@maths.ox.ac.uk}

\author[A.~Swann]{Andrew Swann}
\address[Swann]{Department of Mathematics\\
Aarhus University\\
Ny Munkegade 118, Bldg 1530\\
DK-8000 Aarhus C\\
Denmark\\
\textit{and}\\
CP\textsuperscript3-Origins,
Centre of Excellence for Cosmology and Particle Physics Phenomenology\\
University of Southern Denmark\\
Campusvej 55\\
DK-5230 Odense M\\
Denmark} \email{swann@math.au.dk}

\subjclass[2000]{53C26, 53D20, 14L24}

\begin{abstract}
  We review the quiver descriptions of symplectic and hyperk\"ahler
  implosion in the case of SU(n) actions. We give quiver descriptions
  of symplectic implosion for other classical groups, and discuss some
  of the issues involved in obtaining a similar description for
  hyperk\"ahler implosion.
\end{abstract}

\maketitle

\section{Introduction}

Symplectic implosion is an abelianisation construction in symplectic
geometry invented by Guillemin, Jeffrey and Sjamaar
\cite{Guillemin-JS:implosion}.  Given a symplectic manifold \( M \)
with a Hamiltonian action of a compact group \( K \), its imploded
cross-section \( M_\impl \) is a symplectic stratified space with a
Hamiltonian action of a maximal torus \( T \) of \( K \), such that
the symplectic reductions of \( M \) by \( K \) agree with the
symplectic reductions of the implosion by \( T \).

There is a universal example of symplectic implosion, obtained by
taking \( M \) to be the cotangent bundle \( T^*K \). The imploded
space \( (T^*K)_\impl \) carries a Hamiltonian torus action for which
the symplectic reductions are the coadjoint orbits of \( K \).  It
also carries a Hamiltonian action of \( K \) which commutes with the
\( T \) action, and the implosion \( M_\impl \) of any symplectic
manifold \( M \) with a Hamiltonian action of \( K \) can be
constructed as the symplectic reduction at 0 of the product \( M
\times (T^*K)_\impl \) by the diagonal action of \( K \).

The universal symplectic implosion \( (T^*K)_\impl \) can also be
described in a more algebraic way, as the geometric invariant theory
quotient \( K_\C /\!/N \) of the complexification \( K_\C \) of \( K
\) by a maximal unipotent subgroup~\( N \). This is the affine variety
\( \Spec( \mathcal{O}(K_\C)^N) \) associated to the algebra of
\( N \)-invariant regular functions on \( K_\C \), and may also be
described as the canonical affine completion of the orbit space \(
K_\C/N \) which is a dense open subset of \( K_\C /\!/N \).

Many constructions in symplectic geometry involving the geometry of
moment maps have analogues in hyperk\"ahler geometry. We recall here
that a hyperk\"ahler structure is given by a Riemannian metric \( g \)
and a triple of complex structures satisfying the quaternionic
relations. In fact we then acquire a whole two-sphere's worth of
complex structures, parametrised by the unit sphere in the imaginary
quaternions.  The metric is required to be K\"ahler with respect to
each of the complex structures. In this way a hyperk\"ahler structure
defines a two-sphere of symplectic structures.

Just as the cotangent bundle \( T^*K \) of a compact Lie group carries
a natural symplectic structure, so, by work of Kronheimer, the
cotangent bundle \( T^*K_{\C} \) of the complexified group carries a
hyperk\"ahler structure \cite{Kronheimer:cotangent}. Moreover, in a
series of papers Kronheimer, Biquard and Kovalev showed that the
coadjoint orbits of \( K_{\C} \) admit hyperk\"ahler structures
\cite{Kronheimer:nilpotent,Kronheimer:semi-simple,Biquard,Kovalev}. These
orbits are not however closed in \( \kf_\C^* \) (and the hyperk\"ahler
metrics are not complete) except in the case of semisimple orbits.

In \cite{DKS} and subsequent papers \cite{DKS-Seshadri,DKS-twistor} we
developed a notion of a universal hyperk\"ahler implosion \( (T^*K_\C
)_\hkimpl \) for \( \SU(n) \) actions.  The hyperk\"ahler implosion of
a general hyperk\"ahler manifold \( M \) with a Hamiltonian action of
\( K = \SU(n) \) can then be defined as the hyperk\"ahler quotient of
\( M \times (T^*K_\C )_\hkimpl \) by the diagonal action of \( K \).
As in the symplectic case the universal hyperk\"ahler implosion
carries an action of \( K \times T \) where \( K=\SU(n) \) and \( T \)
is its standard maximal torus. As coadjoint orbits for the complex
group are no longer closed in general, and are not uniquely determined
by eigenvalues, the hyperk\"ahler quotients of \( (T^*K_\C )_\hkimpl
\) by the torus action need not be single orbits. Instead, they are
the Kostant varieties, that is, the varieties in \( \sln(n, \C)^* \)
obtained by fixing the values of the invariant polynomials for this
Lie algebra.  These varieties are unions of coadjoint orbits and are
closures in \( \sln(n,\C)^* \) of the regular coadjoint orbits of \(
K_\C = \SL(n,\C) \). We refer to \cite{Chriss-G:representation},
\cite{Kostant:polynomial} for more background on the Kostant
varieties.

Again by analogy with the symplectic case, we 
can describe the hyperk\"ahler implosion in terms of 
geometric invariant theory (GIT) quotients
by nonreductive group actions. Explicitly, the implosion
is \( (\SL(n ,\C) \times \n^0) \symp N \) where \( N \) is a 
maximal unipotent subgroup of \( K_\C = \SL(n,\C) \) and 
\( \n^0 \) is the annihilator in \( \sln(n,\C)^* \) of its Lie algebra
\( \n \). 
Thus the universal hyperk\"ahler implosion for \( K=\SU(n) \) can be
identified with the complex-symplectic quotient 
\( (\SL(n, \C) \times \n^\circ) \symp N \) of \( T^* \SL(n,\C) \)
by \( N \) in the GIT sense, just as the symplectic implosion is
the GIT quotient of \( K_\C \) by \( N \).

In the case of \( K=\SU(n) \) dealt with in \cite{DKS}, it is possible
to describe the hyperk\"ahler implosion via a purely
finite-dimensional construction using quiver diagrams.  This
construction was motivated by a quiver description of the symplectic
implosion for \( \SU(n) \) we described in \S4 of~\cite{DKS}.

In this article we shall extend these results concerning symplectic
implosion to other classical groups, that is the special orthogonal
and symplectic groups.  Our approach will be inspired by the
description by Lian and Yau \cite{LY} of coadjoint orbits for compact
classical groups using quivers. This suggests ways to extend the
quiver construction of the universal hyperk\"ahler implosion from the
case of \( \SU(n) \) to general classical groups, and we discuss some
of the issues involved here.

\subsubsection*{Acknowledgements.} 
We thank Kevin McGerty and Tom Nevins for bringing the
paper~\cite{GinzburgRiche} to our attention.  The second author is
partially supported by  Swiss National Science Foundation 
grant 200021 138071.  The
fourth author is partially supported by the Danish Council for
Independent Research, Natural Sciences.

\section{Symplectic quivers}
\label{sec:symplectic-quivers}

We begin by trying to construct a quiver model for the universal
symplectic implosion in the case of the orthogonal and symplectic
groups, as was done in \cite{DKS} for special unitary groups.

We consider diagrams of vector
spaces and linear maps
\begin{equation}
  \label{eq:symplectic}
  0 = V_0 \stackrel{\alpha_0}{\rightarrow}
  V_1 \stackrel{\alpha_1}{\rightarrow}
  V_2 \stackrel{\alpha_2}{\rightarrow} \dots
  \stackrel{\alpha_{r-2}}{\rightarrow} V_{r-1}
  \stackrel{\alpha_{r-1}}{\rightarrow} 
  V_{r}= \C^n.
\end{equation}
The dimension vector is defined to be \(
\mathbf n = (n_1,\dots,n_{r-1} ,n_r=n) \) where \( n_i = \dim V_i \).
We will say that the representation is
\emph{ordered} if \( 0\leqslant n_1 \leqslant n_2 \leqslant \dots
\leqslant n_{r} = n \) and \emph{strictly ordered} if \( 0 < n_1
< n_2 < \dots < n_{r} = n \).

In \cite{DKS} we considered \( V_r = \C^n \) as a representation of \(
\SU(n) \), or its complexification \( \SL(n, \C) \). In this setting we
say the quiver is \emph{full flag} if \( r=n \) and \( n_i = i \) for
each \( i \).  We took the geometric invariant theory quotient of the
space of full flag quivers by \( \SL \coloneqq \prod_{i=1}^{r-1} \SL(V_i) \)
(or equivalently, the symplectic quotient by \( \prod_{i=1}^{r-1}
\SU(V_i) \)). The stability conditions imply that the quiver decomposes
into a quiver with zero maps and a quiver with all maps injective. It
was therefore sufficient to analyse the injective quivers up to
equivalence.  We found that the quotient could be stratified into \(
2^{n-1} \) strata, indexing the flag of dimensions of the injective
quivers (after the quivers had been contracted to remove edges where
the maps were isomorphisms). Equivalently, the strata were indexed by
the ordered partitions of \( n \). Each stratum could be identified
with \( \SL(n,\C)/[P,P] \) where \( P \) is the parabolic associated to
the given flag. The upshot was that the full GIT quotient can be
identified as an affine variety with the affine completion \( \SL(n,\C)
\symp N = \Spec \cO(\SL(n,\C))^N \) of the open stratum \(
\SL(n, \C)/N \).

We now wish to
view \( V_r \) as a representation of an orthogonal or symplectic group.
This involves introducing the associated bilinear forms. Our
approach will be motivated by the 
description due to Lian and Yau in \cite{LY} of a quiver approach to generalised
flag varieties for symplectic and orthogonal groups.

Note that for consistency with \cite{DKS} we have altered the notation
of \cite{LY} in some respects. In particular we use \( r \) rather than
\( r+1 \) for the top index, and we use \( n_i \) rather than \( d_i \) for the
dimensions.
 
In the orthogonal case, we let \( J \) denote the matrix with entries
\begin{equation*}
  J_{ij} = \delta_{n+1-i,j} \qquad (1 \leqslant i,j \leqslant n)
\end{equation*}
which are \( 1 \) on the antidiagonal and \( 0 \) elsewhere.  We
therefore have on \( \C^n \) a symmetric bilinear form \( B(v,w) = v^t
Jw \) where \(v^t\) denotes the transpose of \( v \), which is preserved by
\begin{equation*}
  \SO(n,\C) = \{\, g : g^t J g = J \,\}
\end{equation*}
We note that the condition for \( h \) to be in the Lie algebra \( \so(n, \C) \)
is \( h^t J + Jh =0 \), that is, that \( h \) is skew-symmetric about the
ANTI-diagonal. In particular, \( h \) may have arbitrary elements 
in the top left \( d \times d \) block as long as \( d \leqslant \frac{n}{2} \).

Motivated by \cite{LY} let us now consider ordered diagrams where 
\begin{equation*}
n_{r-1} \leqslant \frac{n}{2}
\end{equation*}
and we impose on \( \alpha_{r-1} \) the condition
\begin{equation} \label{cquad}
\alpha_{r-1}^t J \alpha_{r-1} =0.
\end{equation}
Equivalently, this is the condition that the image of \( \alpha_{r-1} \) 
be an isotropic
subspace of \( \C^n \) with respect to \( J \) (which is the reason for the inequality
above).
The space of \( \alpha_{r-1} \) satisfying this condition is
of dimension \( n  n_{r-1} - \frac{1}{2}n_{r-1} (n_{r-1} +1) \).

We let
\(  R(\mathbf n) \) 
be the space of all such diagrams 
satisfying \eqref{cquad} with dimension vector \( \mathbf n \).

Observe that the complexification
 \( \GL \coloneqq \prod_{i=1}^{r-1} \GL(V_i) \)
of 
 \( \tilde{H} \coloneqq \prod_{i=1}^{r-1} \Un(V_i) \)
acts on \( R(\mathbf n) \) by
\begin{align*}
  \alpha_i &\mapsto g_{i+1} \alpha_i g_i^{-1} \quad (i = 1,\dots, r-2),\\
  \alpha_{r-1} &\mapsto \alpha_{r-1} g_{r-1}^{-1}.
\end{align*}
There is also a commuting action of \( \SO(n,\C) \) by left
multiplication of \( \alpha_{r-1} \); note that the full group \(
\GL(n, \C) \) does not now act because it does not preserve
\eqref{cquad}.

As in \cite{DKS}, we shall study the
symplectic quotient of \( R(\mathbf n) \) by the
action of
\begin{equation*}
  H \coloneqq \prod_{i=1}^{r-1} \SU(V_i)
\end{equation*}
or equivalently the GIT quotient of \( R(\mathbf n) \) by its
complexification 
\begin{equation*}
  H_\C = \SL \coloneqq \prod_{i=1}^{r-1} \SL(V_i),
\end{equation*}
viewed as a subgroup of \( \GL \) in the obvious way. This quotient
will have residual actions of the \( r-1 \)-dimensional compact torus
\( T^{r-1} = (S^1)^{r-1} \) and its complexification, as well as of \(
\SO(n, \C) \).

Let us observe that the dimension \( n n_{r-1} - \frac{1}{2}n_{r-1}
(n_{r-1} +1) \) of the set of \( \alpha_{r-1} \) satisfying
\eqref{cquad} equals the dimension of the coset space \( \SO(n,
\C)/\SO(n-n_{r-1}, \C) \). In fact, in the orthogonal case with \( n
\) odd, we can show that this coset space equals the set of injective
\( \alpha_{r-1} \) satisfying~\eqref{cquad}.

For if \( n \) is odd \( \SO(n, \C) \) acts transitively on the set of 
isotropic subspaces of \( \C^n \) of fixed dimension, so \( \alpha_{r-1} \)
can be put into the form
\begin{equation*}
  \begin{pmatrix}
    A_{n_{r-1} \times n_{r-1}} \\
    0_{(n-n_{r-1}) \times n_{r-1}}
  \end{pmatrix}
\end{equation*}
via the \( \SO(n, \C) \) action. As \( n_{r-1} \leqslant \frac{n}{2} \), we
can consider matrices in \( \SO(n, \C) \) with an arbitrary invertible
\( n_{r-1} \times n_{r-1} \) block in the top left corner and a zero
\( (n-n_{r-1}) \times n_{r-1} \) block in the lower left. So in fact
\( \alpha_{r-1} \) can be put into the standard form (used also in the
\( A_n \) case in \cite{DKS})
\begin{equation*}
  \begin{pmatrix}
    I_{n_{r-1} \times n_{r-1}} \\
    0_{(n-n_{r-1}) \times n_{r-1}}
  \end{pmatrix}.
\end{equation*}
The connected component of the stabiliser of this configuration for
the \( \SO(n, \C) \) action is \( \SO(n-n_{r-1}, \C) \).

We now obtain a description of quivers in \( R (\mathbf n) \) with all
\( \alpha_i \) injective, modulo the action of \( \SL \). For,
combining the above observation with the arguments of \S4
of~\cite{DKS}, the action of \( \SL \times \SO(n, \C) \) can be used to
put the maps in standard form
\begin{equation*}
  \alpha_i =
  \begin{pmatrix}
    I_{n_i \times n_i} \\
    0_{(n_{i+1}-n_i) \times n_i}
  \end{pmatrix}.
\end{equation*}
The remaining freedom is a commutator of a parabolic in 
\( \SO(n, \C) \) where the first \( r-1 \) block sizes are \( n_{i+1} - n_i \).
The blocks in the Levi subgroup lie in \( \SL(n_{i+1}-n_i) \)
which is why we get the commutator rather than the full parabolic.
Hence the injective quivers with fixed dimension vector 
\( \mathbf n \) modulo 
the action of \( \SL \) are parametrised by \( \SO(n, \C)/[P,P] \), where
\( P \) is the parabolic associated to the dimension vector.

Note that the blocks corresponding to the upper left square of size
\( n_{r-1} =\sum_{i=0}^{r-2} n_{i+1} - n_i \) will determine the 
blocks in the lower right square of size \( n_{r-1} \), 
at least on Lie algebra level, by the property of being in the orthogonal 
group. 

\begin{remark}
  By intersecting the parabolic (rather than its commutator) with the
  compact group \( \SO(n) \) we get \( \prod_{i=0}^{r-2} \Un(n_{i+1} -
  n_i) \times \SO(n - n_{r-1}) \), which is the isotropy group for the
  associated compact flag variety.  Putting \( p_{i+1} = n_{i+1} - n_i
  \), and \( \ell = n-n_{r-1} \), we get \( \sum_{i=1}^{r-1} p_i = n-l
  \), in accordance with the results of \cite[p.~233, section 8H]{Besse}.
\end{remark}

\medskip For a model for the non-reductive GIT quotient \( K_\C /\!/N
\) in the \( B_k \) case, that is when \( n=2k+1 \) and \( K=\SO(2k+1)
\), we can take \( \mathbf n \) equal to \( (1, 2, 3, \dots, k, 2k+1)
\) which will be the full flag condition in this context. We now
consider the GIT quotient \( R(\mathbf n) \symp \SL \). Note that at
this stage the maps \( \alpha_i \) are not assumed to be injective.

\medskip
As the \( \SL \) action is the same as in the \( A_n \) case,
the stability analysis proceeds as in \cite{DKS}. 
We find that for polystable configurations we may decompose
each vector space \( \C^i \) as
\begin{equation}
  \label{symp:decomp}
  \C^i = \ker \alpha_i \oplus \C^{m_i},
\end{equation}
where \( \C^{m_i} = \im \alpha_{i-1} \) if \( m_i \neq 0 \).
So the quiver decomposes into a zero quiver and an injective quiver.
After contracting legs of the quiver which are isomorphisms, as 
in \cite{DKS}, we obtain a strictly ordered injective quiver
of the form considered above. We stratify the quotient
\( R(\mathbf n) \symp \SL \) by the flag of dimensions of the
injective quiver, as in the \( \SL(n, \C) \) case.

\medskip
We have thus identified the strata of the GIT quotient of the space of
full flag quivers with the strata \( \SO(n, \C)/[P,P] \)
of the universal symplectic implosion, or equivalently
 the non-reductive GIT quotient \( K_\C /\!/N \)  (where \( N \) is a maximal
unipotent subgroup). As the complement
of the open stratum \( \SO(n,\C)/N \) is of complex codimension strictly greater than one, we see
that the implosion \( K_\C /\!/N \) and the GIT quotient of the space
of full flag quivers are affine varieties with the same coordinate
ring \( \cO(\SO(n,\C))^N \), and so are isomorphic.

Guillemin, Jeffrey and Sjamaar \cite{Guillemin-JS:implosion} showed that
the non-reductive GIT quotient \( K_\C /\!/N \) has a \( K \times T \)-invariant
K\"ahler structure such that it can be identified symplectically with
the universal symplectic implosion for \( K \). In order to see this K\"ahler
structure on \( R(\mathbf n) /\!/ \SL \) we can put an \( \tilde{H} \times K \)-invariant
flat K\"ahler structure on \( R(\mathbf n) \) and identify the GIT quotient
\( R(\mathbf n) /\!/ \SL \) with the symplectic quotient \( R(\mathbf n) \symp H \).
To achieve \( \tilde{H} \times K \)-invariance we use the standard
flat K\"ahler structure on \( (\C^{j-1})^* \otimes \C^j \) for \( j \leqslant r-1 \)
but the flat K\"ahler structure defined by \( J \) on 
\begin{equation*}
(\C^{r-1})^* \otimes \C^n \cong (\C^n)^{r-1}.
\end{equation*}

Recall that a polystable quiver decomposes into the sum of a zero quiver and an 
injective quiver, and determines for us a partial flag in \( \C^n \)
\begin{equation*}
W_1 \subseteq W_2 \subseteq \dots \subseteq W_{r-1} \subseteq \C^n,
\end{equation*}
where \( W_j = \im \alpha_{r-1} \circ \alpha_{r-2} \circ
\dots \circ \alpha_j \), whose dimension vector \( (w_1 = \dim W_1,
\dots, w_{r-1} = \dim W_{r-1}) \) is determined by the injective
summand. The condition that \( \alpha_{r-1}^t J \alpha_{r-1} = 0 \)
ensures that \( W_{r-1} \) is an isotropic subspace of \( \C^n \), so
we can use the action of the compact group \( K=\SO(n) \) to put this
flag into standard form with \( W_j \) spanned by the first \( w_j \)
vectors in the standard basis for \( \C^n \).  Then we can use the
action of \( \SL \) to put the polystable quiver into the form where
\begin{equation} \label{quiverT}
  \alpha_j =
  \begin{pmatrix}
    0 & 0 & \dots &  0 \\
    \nu_1^j & 0 & \dots &  0 \\
    0 & \nu_2^j & \dots & 0 \\
    \vdots & \vdots & \ddots & \vdots \\
    0 & 0 &\dots &  \nu^j_j
  \end{pmatrix}
\end{equation}
for \( j<r-1 \) and 
\begin{equation*}
\alpha_{r-1} = 
  \begin{pmatrix}
    0 & 0 & \dots & 0 \\
    \vdots & \vdots & & \vdots \\
    0 & 0 & \dots &  0 \\
    \nu_1^{r-1} & 0 & \dots &  0 \\
    0 & \nu_2^{r-1} & \dots & 0 \\
    \vdots & \vdots & \ddots & \vdots \\
    0 & 0 &\dots &  \nu^{r-1}_{r-1}
  \end{pmatrix}
\end{equation*}
with \( \nu_i^j \in \C \). Let \( R_T(\mathbf n) \) denote the subspace of
\( R(\mathbf n) \) consisting of quivers of this form.

Note that the moment map for the action of the unitary group \( \Un(V_j) \) on
\( R(\mathbf n) \) takes a quiver \eqref{eq:symplectic} to
\begin{equation*}
\bar{\alpha}_j^t \alpha_j - \alpha_{j-1} \bar{\alpha}_{j-1}^t
\end{equation*}
for \( 1 \leqslant j \leqslant r-1 \), so the moment map for the action of the
product \( \tilde{H} = \prod_{j=1}^{r-1} \Un(V_j) \) takes \( R_T(\mathbf n) \)
into the Lie algebra of the product of the standard (diagonal) maximal 
tori \( T_{V_j} \) of the unitary groups \( \Un(V_j) \). Thus there is a natural
map of symplectic quotients
\begin{equation*}
  \theta_T\colon R_T(\mathbf n) \symp H_T \to R(\mathbf n) \symp H
\end{equation*}
where \( H_T = \prod_{j=1}^{r-1} (T_{V_j} \cap \SL(V_j))  \) is a maximal
torus of \( H \). Moreover
\begin{equation*}
R(\mathbf n)\symp H = K \, \theta_T( R_T(\mathbf n) \symp H_T)
\end{equation*}
where \( K = \SO(n) \) and \( R_T(\mathbf n)\symp H_T \) is a toric variety.

The moment map for the action of the torus \( T^{r-1} = (S^1)^{r-1} \) on \( R(\mathbf n)\symp H \)
takes a point  represented by a quiver of the form \eqref{quiverT} satisfying the
moment map equations
\begin{equation*}
  \begin{pmatrix}
    \abs{\nu_1^j}^2 & 0 & \dots &   0 \\
    0 & \abs{\nu_2^j}^2 & \dots  & 0 \\
    \vdots & \vdots & \ddots  & \vdots\\
    0 & 0 & \dots & \abs{\nu^j_j}^2
  \end{pmatrix} 
  = \lambda^\R_j I +
  \begin{pmatrix}
    0 & 0 & \dots &  0 \\
    0 & \abs{\nu_1^{j-1}}^2  & \dots &   0 \\
    \vdots & \vdots & \ddots & \vdots  \\
    0 & 0 & \dots & \abs{\nu^{j-1}_{j-1}}^2
  \end{pmatrix}
\end{equation*}
for some \( \lambda_1^\R, \dots, \lambda_{r-1}^\R \in \R \),
or equivalently  
\begin{equation}
  \label{sympmmap}
  \abs{\nu_i^j}^2 = \lambda_j^\R +
  \lambda^\R_{j-1} + \cdots + \lambda^\R_{j-i+1} \qquad\text{if}\ 1
  \leqslant i \leqslant j < n,
\end{equation}
to \( ( \lambda_1^\R, \dots , \lambda_{r-1}^\R) \) in the Lie algebra
of \( T^{r-1} \), while the moment map for the action of \( K \) takes
this point to
\begin{equation*}
  \begin{pmatrix}
    -\abs{\nu_{r-1}^{r-1}}^2 & \dots & 0 & 0 & \dots & 0 \\
    \vdots & \ddots & \vdots & \vdots & & \vdots \\
    0 & \dots & -\abs{\nu_{1}^{r-1}}^2 & 0 & \dots   & 0 \\
    0 & \dots & 0 &\abs{\nu_{1}^{r-1}}^2 &  \dots &   0 \\
    \vdots & & \vdots & \vdots & \ddots & \vdots \\
    0 & \dots & 0 & 0 & \dots &\abs{\nu^{r-1}_{r-1}}^2
  \end{pmatrix},
\end{equation*}
up to constant scalar factors depending on conventions. The image of the
toric variety \( R_T(\mathbf n)\symp H_T \) under this moment map is the
positive Weyl chamber \( \lie t_+ \) of \( K = \SO(n) \), and we obtain a symplectic 
identification of \( R(\mathbf n)\symp H \) with the universal symplectic
implosion
\begin{equation*}
(T^*K)_{\impl} = (K \times \lie t_+) / \sim
\end{equation*}
of \( K = \SO(n) \), where \( (k,\xi) \sim ( k', \xi') \) if and only if \( \xi = \xi' \),
with stabiliser \( K_\xi \) under the coadjoint action of \( K \), and \( k = k' \tilde{k} \)
for some \( \tilde{k} \in [ K_\xi, K_\xi] \).

\bigskip
We may argue in a very similar way for the symplectic group
\( \SP(2k, \C) \), the complexification of \( \SP(2k) \). Following \cite{LY} we replace the symmetric bilinear form
by the skew form on \( \C^n = \C^{2k} \) defined by the matrix
\begin{equation*}
  J_2 = \begin{pmatrix}
    0  & J \\
    -J & 0
  \end{pmatrix}
\end{equation*}
Once again we find that \( \SP(2k, \C) \) acts transitively on the set
of \( \alpha_{r-1} \) satisfying the condition
\begin{equation*}
\alpha_{r-1}^t J_2 \alpha_{r-1} =0
\end{equation*}
and the above arguments go through mutatis mutandis.

For \( \SO(n, \C) \) with \( n \) even, the isotropic subspace \( \im
\alpha_r \) may be self-dual or anti-self-dual if \( n_{r-1} =
\frac{n}{2} \). We take the component of the locus defined by
\eqref{cquad} corresponding to the image being self-dual, and now we
get the desired transitivity.

\begin{theorem}
  \label{thm:sympimp}
  Let \( K \) be a compact classical group and let us consider full
  flag quivers for \( K_{\C} \) as above. That is, we take \( {\bf n}=
  (n_1, \dots, n_r) \) to be \( (1,2, \dots, n) \) for \( K=\SU(n) \),
  \( (1,2, \dots, k,2k+1) \) for \( \SO(2k+1) \), and \( (1,2,\dots,
  k,2k) \) for \( \SO(2k) \) or \( \SP(k) \). Also in the orthogonal
  and symplectic cases we impose the appropriate isotropy condition on
  the top map in the quiver, and take the appropriate component of the
  space of isotropic subspaces in the even orthogonal case, to obtain
  a space \( R(\mathbf n) \) of full flag quivers.

  Then the symplectic quotient of \( R(\mathbf n) \) by \( H(\mathbf n) =
  \prod_{i=2}^{r-1} \SU(n_i,\C) \) can be identified naturally with the 
universal symplectic implosion for \( K \), or equivalently with the non-reductive
GIT quotient \( K_\C /\!/N \).
  The stratification by quiver diagrams as above
  corresponds to the stratification of the universal symplectic implosion as the disjoint
  union over the standard parabolic subgroups \( P \) of \( K_{\C}
  \) of the varieties \( K_{\C}/[P,P] \).
\end{theorem}

\begin{example} \label{SO3symp}
The lowest rank case of the above construction is when the group is
\( \SO(3) \). The quiver is now just
\begin{equation*}
0 \stackrel{\alpha_0}{\rightarrow} 
\C \stackrel{\alpha_1}{\rightarrow} \C^3
\end{equation*}
where \( \alpha_0 =0 \). As \( H= \SU(1) \) here there is
no quotienting to perform. The matrix \( J \) is
\( \left(\begin{smallmatrix}
0 & 0 & 1 \\
0 & 1 & 0 \\
1 & 0 & 0
\end{smallmatrix}
\right) \) and, putting
\( \alpha_1 = (x,y,z) \), the isotropy condition
\( \alpha_1^t J \alpha_1 =0 \) becomes
\begin{equation*}
y^2 + xz =0.
\end{equation*}
This affine surface is a well known description of the
Kleinian singularity \( \C^2 /\Z_2 \). This is a
valid description of the symplectic implosion for
\( \SO(3) \), since the implosion for the double cover \( \SU(2) \)
is just \( \C^2 \).
\end{example}

\begin{remark}
  We mention here an alternative description of symplectic implosions
  using the concept of Cox rings \cite{COX,HuKeel,LafV}.  If \( X \)
  is an algebraic variety and \( L_1, \dots, L_n \) are generators
  for Pic\( (X) \), then we form the Cox ring
  \begin{equation} \label{Cox} {\rm Cox}(X, L) = \bigoplus_{(m_1,
    \dots, m_n) \in \Z^n} H^0 ( X, m_1 L_1 + \dots + m_n L_n)
  \end{equation}
  Hu and Keel \cite{HuKeel} introduced the class of \emph{Mori dream
  spaces} -- the varieties \( X \) whose Cox ring is finitely
  generated. These include toric varieties, which are characterised by
  Cox\( (X) \) being a polynomial ring. It was proved in \cite{HuKeel}
  that, as the name suggests, Mori dream spaces are well behaved from
  the point of view of the Minimal Model Programme. After a finite
  sequence of flips and divisorial contractions we arrive at a space
  birational to \( X \) which either is a Mori fibre space or has nef
  canonical divisor.  Mori dream spaces may be realised as GIT
  quotients by tori of the affine varieties associated to their Cox
  rings.  The above sequence of flips and contractions can be
  expressed in terms of explicit variation of GIT wall-crossings, and
  indeed the Mori chambers admit a natural identification with
  variation of GIT chambers. Since torus variation of GIT is
  well-understood, in principle the Mori theory of a Mori dream space
  is also well-understood, at least given an explicit enough
  presentation of Cox(X).

  If \( K \) is a compact Lie group and \( P \) is a parabolic
  subgroup of \( K_\C \), then the Cox ring of \( K_{\C}/P \) is the
  coordinate ring of the quasi-affine variety \( K_{\C}/[P,P] \),
  which is finitely generated. In particular \( K_{\C}/P \) is a Mori
  dream space.

  Taking \( P \) to be a Borel subgroup \( B \), we find that the Cox
  ring of \( K_{\C}/B \) is the finitely generated ring \(
  \cO(K_{\C}/N) = \cO(K_\C)^N \) whose associated affine variety \(
  \Spec(\mathcal{O}(K_\C)^N) \) is the universal symplectic implosion
  \( K_{\C} /\!/ N \).
\end{remark}

\section{Hyperk\"ahler quiver diagrams}
\label{sec:hyperk-quiv-diagr}

For \( K=\SU(n) \) actions we developed in \cite{DKS} 
a finite-dimensional approach to
constructing the universal hyperk\"ahler implosion for \( K \) via quiver diagrams.
In that case the symplectic quivers formed a linear space and
we just took the cotangent bundle, which amounted to
putting in maps \( \beta_i \colon V_{i+1} \rightarrow V_i \) in addition to the
\( \alpha_i \colon V_i \rightarrow V_{i+1} \). Writing \( V_i = \C^{n_i} \),
we thus worked with the
flat hyperk\"ahler space
\begin{equation}
  \label{eq:Mn}
  M = M(\mathbf n) = \bigoplus_{i=1}^{r-1} \HH^{n_i n_{i+1}} = 
  \bigoplus_{i=1}^{r-1} \Hom(\C^{n_i},\C^{n_{i+1}}) \oplus
  \Hom(\C^{n_{i+1}},\C^{n_i})
\end{equation}
with the hyperk\"ahler action of \( \Un(n_1) \times \dots \times
\Un(n_r) \)
\begin{equation*}
  \alpha_i \mapsto g_{i+1} \alpha_i g_i^{-1},\quad
  \beta_i \mapsto g_i \beta_i g_{i+1}^{-1} \qquad (i=1,\dots r-1),
\end{equation*}
with \( g_i \in \Un(n_i) \) for \( i=1, \dots, r \). 
Right quaternion multiplication was given by
\begin{equation}
  \label{eq:j}
  (\alpha_i,\beta_i)\jj = (-\beta_i^*,\alpha_i^*).
\end{equation}
If each \( \beta_i \) is zero we recovered a
symplectic quiver diagram.

We considered the hyperk\"ahler quotient of \( M(\mathbf n) \) with respect to
the group \( H = \prod_{i=1}^{r-1}\SU(n_i) \), obtaining a stratified
hyperk\"ahler space \( Q = M \hkq H \), which has a residual action of
the torus \( T^{r-1} = \tH / H \) where \( \tH = \prod_{i=1}^{r-1}\Un(n_i) \), 
as well as a commuting action of \( \SU(n_r) = \SU(n) \). 
When \( \mathbf{n} = (1,2,\dots, n) \) we can identify this torus with 
the standard maximal torus \( T \) of \( \SU(n) \) using the simple roots of \( T \).

  The \emph{universal hyperk\"ahler implosion for \( \SU(n) \)} is defined
to  be the hyperk\"ahler quotient \( Q = M \hkq H \), where \( M \), \(
  H \) are as above with \( n_j = j \), for \( j=1, \dots, n
  \), (i.e.\ the case of a full flag quiver).

From the complex-symplectic viewpoint, \( Q \)
is  the GIT quotient, by the 
complexification
\begin{equation*}
  H_\C  = \prod_{i=1}^{r-1}\SL(n_i,\C)
\end{equation*}
of \( H \), of the zero locus of the complex moment map \( \mu_{\C} \) for the
\( H \) action. 

The components of this complex moment map \( \mu_\C \) are given by
the tracefree parts of \( \alpha_{i-1}\beta_{i-1} - \beta_i \alpha_i
\).
The complex moment map equation \( \mu_\C =0 \) can thus be expressed
as saying
\begin{equation}
  \label{eq:mmcomplex}
   \beta_i \alpha_i - \alpha_{i-1}\beta_{i-1}  = \lambda^\C_i I \qquad
  (i=1,\dots,r-1),
\end{equation}
for some complex scalars \( \lambda^\C_1,\dots,\lambda^\C_{r-1} \),
while the real moment map equation is given by
\begin{equation}
  \label{eq:mmreal}
\beta_{i-1}^* \beta_{i-1} -\alpha_{i-1} \alpha_{i-1}^*
 - \beta_i \beta_i^* + \alpha_{i}^* \alpha_{i}  = \lambda^\R_i I \qquad
  (i=1,\dots,r-1),
\end{equation}
for some real scalars \( \lambda^\R_1,\dots,\lambda^\R_{r-1} \).

 The action of \( H_\C \) is given by
\begin{gather*}
  \alpha_i \mapsto g_{i+1} \alpha_i g_i^{-1}, \quad \beta_i \mapsto
  g_i \beta_i g_{i+1}^{-1} \qquad (i=1,\dots r-2),\\
  \alpha_{r-1} \mapsto \alpha_{r-1} g_{r-1}^{-1}, \quad \beta_{r-1}
  \mapsto g_{r-1} \beta_{r-1},
\end{gather*}
where \( g_i \in \SL(n_i,\C) \).  
The residual action of \( \SL(n,\C) =
\SL(n_r,\C) \) on the quotient \( Q \) is given by
\begin{equation*}
  \alpha_{r-1} \mapsto g_r \alpha_{r-1}, \quad
  \beta_{r-1} \mapsto \beta_{r-1} g_r^{-1}.
\end{equation*}

There is also a residual action of \( \tH_\C/H_\C \) which we can
identify, in the full flag case, with the maximal torus \( T_\C \) of
\( K_\C \). The complex numbers \( \lambda_i \) combine to give the
complex-symplectic moment map for this complex torus action.  We
remark that reduction of \( Q \) by the maximal torus at level~\( 0 \)
recovers the construction of the nilpotent variety
\cite{Kobak-S:finite}, \cite{KP1}. While there is a similar
quiver description of the nilpotent
variety for the classical algebras $\so(n, \C)$ and $\sP(n, \C)$ 
 \cite{Kobak-S:finite}, \cite{KP2}, the construction of an implosion
does not directly generalise, partly because the corresponding groups
$\tilde H$ do not have sufficiently large centers.

Note that as \( Q \) is a hyperk\"ahler reduction by \( H \) at
level~\( 0 \), it also inherits an \( \SU(2) \) action that rotates the
two-sphere of complex structures (see \cite{DKS} for details).

Given a quiver \( (\alpha,\beta) \in M(\mathbf n) \), the composition 
\begin{equation*}
  X = \alpha_{r-1} \beta_{r-1} \in \Hom (\C^n,\C^n).
\end{equation*}
is invariant under the action of \( \tH_\C \) and transforms by
conjugation under the residual \( \SL(n,\C) \) action.  The map \( Q
\rightarrow \sln(n, \C) \) given by sending \( (\alpha, \beta) \) to
the tracefree part of \( X \) is therefore \( T_{\C} \)-invariant and
\( \SL(n,\C) \)-equivariant.

In \cite{DKS} and \cite{DKS-Seshadri} we introduced stratifications of
the implosion \( Q \), one reflecting its hyperk\"ahler structure and
one reflecting the group structure of \( \SU(n) \).  We recall in
particular that the open subset of \( Q \) consisting of quivers with
all \( \beta \) surjective may be identified with \( \SL(n, \C)
\times_N \n^0 \cong \SL(n, \C) \times_N \bmf \).

The open stratum \( Q^{hks} \) in the hyperk\"ahler stratification of \( Q \)
consists of the quivers which are hyperk\"ahler stable; that is, for a generic
choice of complex structure all the maps \( \alpha_i \) are injective and
all the maps \( \beta_i \) are surjective. In this situation the kernels of the compositions
\begin{equation*}
\beta_j \circ \beta_{j+1} \circ \cdots \circ \beta_{n-1}
\end{equation*}
for \( 1 \leqslant j \leqslant n \) form a full flag in \( \C^n \); we can use the action of
\( K = \SU(n) \) (which preserves the hyperk\"ahler structure) to put this
flag into standard position. Next we can use the action of 
\( \SL = H_\C \) to put the maps \( \beta_j \) into the form
\begin{equation}
  \label{betaform}
  \beta_j =
  \begin{pmatrix}
    0 & \mu_1^j & 0 & \dots &0  & 0 \\
    0 & 0 & \mu_2^j &  \dots &0 & 0 \\
    \vdots & \vdots & \vdots & \ddots & \vdots & \vdots \\
    0 & 0 & 0 & \dots &\mu_{j-1}^j & 0 \\
    0 & 0 & 0 &\dots &0 & \mu^j_j
  \end{pmatrix}
\end{equation}
for some \( \mu^j_i \in \C \setminus\{ 0 \} \). Then it follows
from the complex moment map equations \eqref{eq:mmcomplex} that
the maps \( \alpha_j \) have the form
\begin{equation}
  \label{alphaform}
  \alpha_j =
  \begin{pmatrix}
    * & * & \dots &* & * \\
    \nu_1^{j} & * & \dots &*  & * \\
    0 &  \nu_2^j & \dots &* & *\\
    \vdots & \vdots & \ddots &\vdots& \vdots \\
    0 &0 & \dots & \nu^{j-1}_{j} &*\\
    0 &0 & \dots & 0&\nu^{j}_{j}
  \end{pmatrix}
\end{equation}
and that the same equations \eqref{eq:mmcomplex} are satisfied if each
\( \alpha_j \) is replaced with
\begin{equation} \label{alphaTform}  
 \alpha_j^t =  
  \begin{pmatrix}
    0 & 0 & \dots &0 & 0 \\
    \nu_1^{j} & 0 & \dots &0  & 0 \\
    0 &  \nu_2^j & \dots &0 & 0\\
    \vdots & \vdots & \ddots &\vdots& \vdots \\
    0 &0 & \dots & \nu^{j-1}_{j} &0\\
    0 &0 & \dots & 0&\nu^{j}_{j}
  \end{pmatrix}
\end{equation}
For a fixed choice of complex structures let us denote by 
\( \bmf^{(\circ)}_+ \) the subset of \( Q \) represented by all quivers of
the form \eqref{alphaform} and \eqref{betaform} satisfying the
hyperk\"ahler moment map equations with \( \mu_i^j \) and \( \nu_i^j \) nonzero
complex numbers. 
Its \( K \)-sweep \( K \bmf^{(\circ)}_+ \) in \( Q \) is then isomorphic to
\begin{equation*}
K \times_T \bmf^{(\circ)}_+ \cong K_\C \times_B \bmf^{(\circ)}_+
\end{equation*}
and consists of all quivers \eqref{eq:symplectic} in \( Q \) such that for each \( j \)  the map \( \alpha_j \)
is injective and the map \( \beta_j \) is surjective and   \( \C^n \) is 
the direct sum of 
\begin{equation*}
\ker (\beta_j \circ \beta_{j+1} \circ \cdots \circ \beta_{n-1})
\end{equation*}
and \( \im (\alpha_{n-1} \circ \cdots \circ \alpha_{j}) \). It follows
that \( K \bmf^{(\circ)}_+ \) is open in \( Q \), and its sweep \(
\SU(2) K \bmf^{(\circ)}_+ \) under the action of \( \SU(2) \) which
rotates the complex structures (and commutes with the action of \( K
\)) is the open stratum \( Q^{hks} \) of \( Q \).

Associating to a quiver \eqref{eq:symplectic} in \( \bmf^{(\circ)}_+
\) with the maps \( \alpha_j \) and \( \beta_j \) in the form
\eqref{alphaform} and \eqref{betaform} the quiver in which \( \alpha_j
\) is replaced with \( \alpha^T_j \) given by \eqref{alphaTform}
defines a map \( \psi \) from \( \bmf^{(\circ)}_+ \) to the hypertoric
variety \( Q_T \) defined in \cite{DKS-Seshadri}.  This hypertoric
variety is the hyperk\"ahler quotient of the space \( M_T \) of all
quivers of the form \eqref{alphaTform} and \eqref{betaform} by the
action of the maximal torus \( H_T \) of \( H \) with the induced
action of \( \tilde{H}_T/H_T = (S^1)^{n-1} \) which is identified with
\( T \) via the basis of \( \lie t^* \cong \lie t \) corresponding to
the simple roots.  The image of the map \( \psi \) is the open subset
\( Q_T^{(\circ)} \) of \( Q_T \) represented by all quivers of this
form with \( \mu_i^j \) and \( \nu_i^j \) all nonzero.

The restriction to \( \bmf^{(\circ)}_+ \) of the complex moment map
for the action of \( K \) associates to a quiver of the form
\eqref{alphaform}, \eqref{betaform} the upper triangular matrix
\begin{equation*}
\alpha_{n-1} \beta_{n-1} - \mathrm{tr}(\alpha_{n-1} \beta_{n-1}) \frac{I}{n}
\end{equation*}
and thus takes values in \( \bmf = \lie t_\C \oplus \n \). Combining
the map \( \psi \) with the projection to \( \n \) of this complex
moment map gives us isomorphisms
\begin{equation*}
\bmf^{(\circ)}_+ \cong Q^\circ_T \times \n
\end{equation*}
and 
\begin{equation*}
K \bmf^{(\circ)}_+ \cong K \times_T (Q^\circ_T \times \n).
\end{equation*}
Under this identification the complex moment map for \( T \) is given
by the (\( T \)-invariant) complex moment map \( \phi\colon Q_T \to \lie
t_\C^* \) for the action of \( T \) on \( Q_T \), and the complex
moment map for \( K \) is given by
\begin{equation*}
[k,\eta, \zeta] \mapsto \mathrm{Ad}^* k \, (\phi(\eta) + \zeta)
\end{equation*}
for \( k \in K \), \( \eta \in Q^\circ_T \) and \( \zeta \in \n \).

The hyperk\"ahler moment map for \( T \) associates to a quiver
satisfying the hyperk\"ahler moment map equations \eqref{eq:mmcomplex}
and \eqref{eq:mmreal} the element \( (\lambda^\C_1, \lambda_1^\R,
\dots, \lambda^\C_{n-1}, \lambda_{n-1}^\R) \) of \( (\C \oplus
\R)^{n-1} \) identified with \( \lie t^* \otimes (\C \oplus \R)
\cong\lie t^* \otimes \R^3 \) via the basis of simple roots. The image
of its restriction to \( Q^{hks} \) is the open subset of \( \lie t^*
\otimes \R^3 \) defined by \( (\lambda_j^\C, \lambda_j^\R) \neq (0,0)
\) for \( j=1,\dots,n-1 \), while the image of \( K \bmf^{(\circ)}_+
\) is the open subset \( (\lie t^* \otimes \R^3)^\circ \) defined by
\( \lambda_j^\C \neq 0 \) for \( j=1, \dots,n-1 \).  Using the same
basis the hypertoric variety \( Q_T \) can be identified with \(
\HH^{n-1} \) and \( Q^{(\circ)}_T \) then corresponds to the
open subset
\begin{equation*}
  \{\,(a_1 + jb_1, \dots, a_{n-1} + j b_{n-1}) \in \HH^{n-1}: a_\ell,
  b_\ell \in \C \setminus \{ 0 \} \,\}. 
\end{equation*}
Under this identification the hyperk\"ahler moment map \( \phi\colon
Q^{(\circ)}_T \to \lie t^* \otimes \R^3 \) is given by
\begin{equation*}
\phi(a_1 + jb_1, \dots, a_{n-1} + j b_{n-1}) =
\end{equation*}
\begin{equation*}
(a_1 b_1, \abs{a_1}^2 - \abs{b_1}^2, \dots, a_{n-1} b_{n-1},
\abs{a_{n-1}}^2 - \abs{b_{n-1}}^2);
\end{equation*}
its fibres are single \( T \)-orbits in \( Q^{(\circ})_T \).

From the description of \( K \bmf^{(\circ)}_+ \) above it follows that the hyperk\"ahler moment
map for \( T \) restricts to a locally trivial fibration
\begin{equation*}
Q^{hks} \to \SU(2) (\lie t^* \otimes \R^3)^\circ
\end{equation*}
over the open subset \( \SU(2) (\lie t^* \otimes \R^3)^\circ  \) of
\( \lie t^* \otimes \R^3  \) with fibre \( K \times \n \).

Similarly the other strata in the hyperk\"ahler stratification of \( Q \) are
constructed from hyperk\"ahler stable quivers of the form 
\begin{equation} \label{Q1}
  0 \stackrel[\beta_0]{\alpha_0}{\rightleftarrows}
  \C^{n_1}\stackrel[\beta_1]{\alpha_1}{\rightleftarrows}
  \C^{n_2}\stackrel[\beta_2]{\alpha_2}{\rightleftarrows}\dots
  \stackrel[\beta_{r-2}]{\alpha_{r-2}}{\rightleftarrows} \C^{n_{r-1}}
  \stackrel[\beta_{r-1}]{\alpha_{r-1}}{\rightleftarrows} \C^{n_r} = \C^n.
\end{equation}
Again for generic choices of complex structures  for each \( j \)  the map \( \alpha_j \)
is injective and the map \( \beta_j \) is surjective and   \( \C^n \) is 
the direct sum of 
\begin{equation*}
\ker (\beta_j \circ \beta_{j+1} \circ \cdots \circ \beta_{r-1})
\end{equation*}
and \( \im (\alpha_{r-1} \circ \cdots \circ \alpha_{j}) \), and we can
use the action of \( K = \SU(n) \) to put the flag in \( \C^n \)
defined by the subspaces \( \ker (\beta_j \circ \beta_{j+1} \circ
\cdots \circ \beta_{r-1}) \) into standard position. Next we can use
the action of \( \prod_{j=1}^{r-1} \SL(n_j) \) to put the maps \(
\beta_j \) into block form of the same shape as \eqref{betaform} where
now each \( \mu^j_i \) is a nonzero scalar multiple of an identity
matrix. Again it follows from the complex moment map equations
\eqref{eq:mmcomplex} that the maps \( \alpha_j \) have block form
similar to \eqref{alphaform} where each \( \nu^j_i \) is a nonzero
scalar multiple of an identity matrix, and that the same equations
\eqref{eq:mmcomplex} are satisfied if each \( \alpha_j \) is replaced
with \( \alpha_j^T \) in block diagonal form as at
\eqref{alphaTform}. We find that the space \( Q_1^{hks} \) of
hyperk\"ahler quivers of the form \eqref{Q1} fibres over an open
subset of \( \lie t_1^* \otimes \R^3 \) (where \( T_1 \cong
(S^1)^{r-1} \)) with fibre
\begin{equation*}
K \times_{[K \cap P_1, K \cap P_1]} \p^0_1
\end{equation*}
where \( P_1 \) is the standard parabolic in \( K_\C \) corresponding
to the flag defined by the subspaces \( \ker (\beta_j \circ
\beta_{j+1} \circ \cdots \circ \beta_{r-1}) \), and \( \p_1^0 \) is
the annihilator in \( \lie k_\C^* \) of its Lie algebra \( \p_1 \).
Note that using the standard pairing on \( \lie k_\C \) and
identifying \( \bmf \) with the annihilator \( \n^0 \) of \( \n \) in
\( \lie k^*_\C \), we have a projection from \( \n \cong \n^* \) onto
the annihilator of \( \p_1 \) in \( \n^* \), and this annihilator can
be identified with \( \p_1^0 \) since \( \n + \p_1 = \lie k_\C \).

By \cite[Proposition 6.9]{DKS} each stratum in \( Q \) can be
identified with a hyperk\"ahler modification
\begin{equation*}
\hat{Q}_1^{hks} = (Q_1^{hks} \times (\HH \setminus \{ 0 \})^\ell) \hkq T^\ell
\end{equation*}
of \( Q_1^{hks} \) for some \( Q_1 \) as above, and the restriction to
this stratum \( \hat{Q}^{hks}_1 \) of the hyperk\"ahler moment map for
\( T \) is a locally trivial fibration over an open subset of
\begin{equation*}
\mathrm{Lie}(Z_K(K \cap P_1))^* \otimes \R^3
\end{equation*}
(where \( Z_K(K \cap P_1) \subseteq T \) is the centre of \( K \cap P_1 \) in \( K \)) with fibre
\begin{equation*}
K \times_{[K \cap P_1, K\cap P_1]} \p_1^0.
\end{equation*}
Using the surjection
\begin{equation} \label{HO} K \times \n \to K \times_{[K \cap P_1, K \cap P_1]} \p_1^0 \end{equation}
induced by the projection \( \n \cong \n^* \to \p_1^0 \), we can lift this locally
trivial fibration to one with fibre \( K \times \n \) which surjects onto the stratum
\( \hat{Q}_1^{hks} \).

In order to patch together these locally trivial fibrations for the
different strata, we can blow up the hypertoric variety \( Q_T \cong
\HH^{n-1} \), replacing it with \( \tilde{Q}_T \cong
\tilde{\HH}^{n-1} \) where \( \tilde{\HH}^{n-1} \) is
the blow-up of \( \HH \cong \C^2 \) at 0 using the complex
structure on \( \HH \) given by right multiplication by \( i
\); this commutes with the hyperk\"ahler complex structures given by
left multiplication by \( i,j \) and \( k \) and also with the action
of the \( S^1 \) component of the maximal torus \( T \cong (S^1)^{n-1}
\).

Note that the hyperk\"ahler moment map for the action of \( S^1 \) on \( \HH \)
induces an identification of the topological quotient \( \HH/S^1 \) with \( \R^3 \);
this pulls back to an identification of \( \tilde{\HH}/S^1 \) with the manifold with
boundary \( \tilde{\R}^3 = (\R^3 \setminus \{ 0 \}) \sqcup S^2 \). Let \( \tilde{Q} \)
be the fibre product
\begin{equation*}
  \begin{CD} \tilde{Q} @>>> Q \\
    @VVV @VVV \\
    \lie t^* \otimes \tilde{R}^3 @>>> \lie t^* \otimes \R^3.
  \end{CD}
\end{equation*}
The descriptions above of the hyperk\"ahler strata of \( Q \) as the images
of surjections from locally trivial fibrations over subsets of
\( \lie t^* \otimes \R^3 \) with fibre \( K \times \n \) patch together to
give a locally trivial fibration
\begin{equation*}
\hat{Q} \to \lie t^* \otimes \tilde{\R}^3
\end{equation*}
with fibre \( K \times \n \) over the manifold with corners \(
\tilde{\R}^3 \), and surjections \( \hat{\chi} \colon \hat{Q} \to
\tilde{Q} \) and \( \chi\colon \tilde{Q} \to Q \) where \( \hat{\chi} \)
collapses fibres via the surjections \eqref{HO} and \( \chi \) is the
pullback of the surjection \( \lie t^* \otimes \tilde{\R}^3 \to \lie
t^* \otimes \R^3 \).

\begin{remark}
If at \eqref{betaform} we only allow ourselves to use the
action of \( H \), not \( H_\C \), to put the maps \( \beta_j \) into standard form,
then we are able to ensure that each \( \beta_j \) is of the form
\begin{equation}
  \label{betaform2}
  \beta_j =
  \begin{pmatrix}
    0      & \mu_1^j & *       & \dots & *           & * \\
    0      & 0       & \mu_2^j & \dots & *           & * \\
    \vdots & \vdots  & \vdots  & \ddots & \vdots  & \vdots \\
    0      & 0       & 0       & \dots & \mu_{j-1}^j & * \\
    0      & 0       & 0       & \dots & 0           & \mu^j_j
  \end{pmatrix}
\end{equation}
for some \( \mu^j_i \in \C \setminus\{ 0 \} \). It still follows
from the complex moment map equations \eqref{eq:mmcomplex} that
the maps \( \alpha_j \) then have the form \eqref{alphaform}.
Similarly an element of any stratum \( \hat{Q}_1^{hks} \) as above
can for generic choices of complex structures be put into
block form \eqref{alphaform} and \eqref{betaform2} using the action of \( K \times H \), where now 
\( \mu_i^j \) and \( \nu_i^j \) denote nonzero scalar multiples of identity matrices.

\end{remark}

\section{Properties of hyperk\"ahler implosion}

In this section we will list some of the main properties of the universal
hyperk\"ahler implosion \( Q = (T^*K_\C)_{\hkimpl} \) for \( K = \SU(n) \) which
we expect to be true for more general compact groups \( K \).

i) \( Q \) is a stratified hyperk\"ahler space of real 
dimension \( 2(\dim K + \dim T) \)
where \( T \) is a maximal torus of \( K \). It has an action of \( K \times T \) which
preserves the hyperk\"ahler structure and has a hyperk\"ahler moment map
\begin{equation*}
  \mu^{K \times T}\colon Q \to (\lie k^* \oplus \lie t^*) \otimes \R^3,
\end{equation*}
as well as a commuting action of \( \SU(2) \) which rotates the complex
structures on \( Q \) (see \cite{DKS} for the case \( K = \SU(n) \)).

ii) The hyperk\"ahler reduction at 0 of \( Q \) by \( T \) can be
identified for any choice of complex structure, via the complex moment
map for the action of \( K \), with the nilpotent cone \( \mathcal{N}
\) in \( \lie k_\C \). We can view this as the statement that the
Springer resolution \( \SL(n, \C) \times_{B} \n \rightarrow {\mathcal
N} \) is an affinisation map.  The reduction at a generic point of \(
\lie t ^*\otimes \R^3 \) is a semisimple coadjoint orbit of \( K_\C
\), and in general the hyperk\"ahler reduction of \( Q \) by \( T \)
at any point of \( \lie t ^*\otimes \R^3 \) can be identified for any
choice of complex structure, via the complex moment map for the action
of \( K \), with a Kostant variety in \( \lie k^*_\C \) (that is, the
closure of a coadjoint orbit).  We refer to \cite{DKS} for the case \(
K = \SU(n) \).

iii) When \( K \) is semisimple, simply connected and connected (as for special unitary
groups) its universal symplectic implosion embeds in the affine space
\begin{equation*}
\bigoplus_{\varpi \in \Pi} V_\varpi,
\end{equation*}
where \( \{\, V_\varpi: \varpi \in \Pi \,\} \) is the set of
fundamental representations of \( K \), as the closure of the \( K_\C
\)-orbit of \( v = \sum_{\varpi \in \Pi} v_\varpi \) for any choice of
highest weight vector \( v_\varpi \) for the irreducible
representation \( V_\varpi \). When \( K = \SU(n) \) it was shown in
\cite{DKS-twistor} that the universal hyperk\"ahler implosion \( Q \)
embeds in the space
\begin{equation*}
H^0(\mathbb{P}^1, ((\lie k^*_\C \oplus \lie t^*_\C) \otimes \cO(2)) \oplus\bigoplus_{\varpi}
 V_{\varpi} \otimes \cO(j(\varpi)))
\end{equation*}
of holomorphic sections of the vector bundle
\begin{equation} \label{vecotb} \mathcal{V} =  
((\lie k^*_\C \oplus \lie t^*_\C) \otimes \cO(2)) \oplus\bigoplus_{\varpi}
 V_{\varpi} \otimes \cO(j(\varpi)) \end{equation}
over \( \mathbb{P}^1 \) for suitable positive integers \( j(\varpi) \). Moreover this embedding induces a
holomorphic and generically injective map from the twistor space \( \mathcal{Z} Q \) of \( Q \) to the
vector bundle \( \mathcal{V} \) over \( \mathbb{P}^1 \), and the hyperk\"ahler structure can be recovered
from this embedding when \( K = \SU(n) \) (\cite{DKS-twistor}).

iv) Let \( N \) be a maximal unipotent subgroup of the complexification
\( K_\C \) of \( K \).
It was shown in \cite{DKS} that when \( K = \SU(n) \) the algebra
of invariants \( \cO(K_\C \times \n^0)^N \)  is finitely generated
and for any choice of complex structures \( Q \) is isomorphic to
the affine variety 
\begin{equation*}
  (K_\C \times \n^0) \symp N =  \Spec \cO (K_\C \times \n^0)^N
\end{equation*}
 associated to this algebra of invariants. This variety may be
viewed as
the complex-symplectic quotient (in the sense of nonreductive GIT) 
of \( T^*K_\C = K_\C \times \lie k_\C^* \) by the action of \( N \)
given by \( (g, \zeta) \mapsto (gn^{-1}, Ad(n) \zeta) \).
With respect to this identification the complex moment maps for the
commuting \( K \) and \( T \) actions on \( Q \) are the morphisms from  
 \( (K_\C \times \n^0) \symp N \) induced by the \( N \)-invariant
morphisms from \( K_\C \times \n^0 \) to \( \lie k^*_\C \) and \( \lie t^*_\C \) given
by
\begin{equation*}
(g,\zeta) \mapsto \mathrm{Ad}^*(g) \zeta
\end{equation*}
and 
\begin{equation*}
(g,\zeta) \mapsto  \zeta_T
\end{equation*}
where \( \zeta_T \in \lie t^*_\C \) is the restriction of \( \zeta \in \n^0 \subseteq \lie k^*_\C \)
to \( \lie t_\C \).

It has been proved very recently by Ginzburg and Riche \cite[Lemma
3.6.2]{GinzburgRiche} that the algebra of regular functions on \(
T^*(G/N) \) is finitely generated for a general reductive \( G \) with
maximal unipotent subgroup \( N \). Taking \( G=K_\C \) for any
compact group \( K \) this cotangent bundle may be identified with \(
K_\C \times_N \n^0 \), and its algebra of regular functions is \(
\cO(K_\C \times \n^0)^N \). Hence the non-reductive GIT quotient \(
(K_\C \times \n^0 ) \symp N = \Spec \cO (K_\C \times \n^0)^N \) is a
well defined affine variety in general, and is the canonical affine
completion of the quasi-affine variety \( T^*(K_\C/N) \) just as \(
K_\C/\!/N \) is the canonical affine completion of the quasi-affine
variety \( K_\C/N \). It is enough to consider the case when \( K \)
is semisimple, connected and simply connected. Then their proof
provides a reasonably explicit set of generators involving the
fundamental representations \( V_\varpi \) of \( K \) and these give
an embedding of the affine variety \( (K_\C \times \n^0 ) \symp N =
\Spec  \cO (K_\C \times \n^0)^N \) as a closed subvariety of the
space of sections \( H^0(\mathbb{P}^1, \mathcal{V}) \) of a vector
bundle \( \mathcal{V} \) over \( \mathbb{P}^1 \) as at \eqref{vecotb}
above.  Note also that the GIT complex-symplectic quotient at level \(
0 \) of \( (K_\C \times \n^0 ) \symp N \) may be viewed as \( (K_\C
\times \n) \symp B \) which is the nilpotent variety (see the remarks
in ii) above).  Similarly reductions at other levels will yield the
Kostant varieties (cf.\ the discussion in \S3 of~\cite{DKS}).

Thus we expect that in general, as in the case when \( K = \SU(n) \),
\( (K_\C \times \n^0 ) \symp N \) has a hyperk\"ahler structure
determined by this embedding and can be identified with the
universal hyperk\"ahler implosion for \( K \).

Note that the scaling action of \( \C^* \) on \( \n^0 \) induces an
action of \( \C^* \) on \( K_\C \times \n^0 \) which commutes with the
action of \( N \) and thus induces an action of \( \C^* \) on \( (K_\C
\times \n^0 ) \symp N \).
Since 
\( \C^* \) acts on \( \cO (K_\C \times \n^0) \), and thus on \( \cO
(K_\C \times \n^0)^N \), with only non-negative weights, the sum of
the strictly positive weight spaces forms an ideal \( I \) in \( \cO
(K_\C \times \n^0)^N \) which defines the fixed point set for the
action of \( \C^* \) on \( (K_\C \times \n^0 ) \symp N \).  This fixed
point set is therefore the affine variety \( \Spec(\cO (K_\C \times
\n^0)^N/I) \), which can be naturally identified with \( \Spec ( (\cO
(K_\C \times \n^0)^N)^{\C^*}) \) and thus with the universal symplectic
implosion \( K_\C \symp N = \Spec (\cO (K_\C)^N) \).
  
v) If \( \zeta = (\zeta_1, \zeta_2, \zeta_3) \in \lie k^* \otimes \R^3
\) let
\begin{equation*}
K_\zeta = K_{\zeta_1} \cap K_{\zeta_2} \cap K_{\zeta_3}
\end{equation*}
where \( K_{\zeta_j} \) is the stabiliser of \( \zeta_j \) under the
coadjoint action, and let \( \mathcal{N}_\zeta \) be the nilpotent
cone in \( (\lie k_\zeta)^*_\C \) which we identify with \( (\lie
k_\zeta)_\C \) as usual. By work of Kronheimer
\cite{Kronheimer:nilpotent} there is a \( K_\zeta \times T \times
\SU(2) \)-equivariant embedding
\begin{equation*}
\mathcal{N}_\zeta \to \lie k_\zeta \otimes \R^3
\end{equation*}
whose composition with the projection from \( \lie k_\zeta \otimes
\R^3 \) to \( (\lie k_\zeta)_\C \) for any choice of complex
structures is the inclusion of the nilpotent cone \( \mathcal{N}_\zeta
\) in \( (\lie k_\zeta)_\C \). From the discussion in \( \S \)3 we
expect that for any compact group \( K \) the image of the
hyperk\"ahler moment map for the action of \( K \) on the universal
hyperk\"ahler implosion \( Q \) should be
the $K$-sweep of
\begin{equation*}
  \lie t_{(\mathrm{hk})} = \{\, \zeta + \xi \in \lie k \otimes \R^3: \zeta \in \lie
  t \otimes \R^3 \ \text{and}\ \xi \in \mathcal{N}_\zeta \,\}
\end{equation*}
and the hyperk\"ahler implosion 
\(  X_{\hkimpl} = (X \times Q) \hkq K \) for any hyperk\"ahler manifold \( X \)
with a Hamiltonian hyperk\"ahler action of \( K \) and hyperk\"ahler moment
map \( \mu_X\colon X \to \lie k \otimes \R^3 \) should be  given by
\begin{equation*}
X_{\hkimpl} = \mu_X^{-1}(\lie t_{(\mathrm{hk})})/\sim.
\end{equation*}
Here \( x \sim y \) if and only if \( \mu_X(x) = \zeta + \xi \) and
\( \mu_X(y) = w(\zeta + \xi') \) for some \( \zeta \in \lie t \otimes \R^3 \), some \( \xi, \xi'
\in \mathcal{N}_\zeta \subseteq \lie k \otimes \R^3 \), some \( w \) in the Weyl group \(W\) of \( K \), identified with a finite subgroup of the normaliser of \( T \) in \( K \), and moreover
\( x = k w^{-1} y \) for some \( k \in [K_\zeta, K_\zeta] \).

\section{Hyperk\"ahler implosion for special orthogonal and symplectic
groups}

In the case of \( K= \SU(n) \) the quiver model for the universal
symplectic implosion is a symplectic quotient of a flat linear space,
so to obtain a quiver model for the universal hyperk\"ahler implosion
we could take its cotangent bundle (replacing symplectic with
hyperk\"ahler quivers) and the corresponding hyperk\"ahler quotient.

We would like to mimic this construction for the orthogonal and symplectic groups.
However we now have the problem that the space of symplectic quivers
has a non-flat piece since the top map \( \alpha_{r-1} \) has to satisfy the
system of quadrics \eqref{cquad} given by
\( \alpha_{r-1}^t J \alpha_{r-1} =0 \).

If this system of equations cut out a smooth variety we could appeal to a result
of Feix \cite{Feix} (see also Kaledin \cite{Kal}) that gives a hyperk\"ahler structure
on an open neighbourhood of the zero section of the cotangent bundle of
a K\"ahler manifold with real-analytic metric. In our case, however, the
variety defined by \eqref{cquad} is singular. We could of course stratify into
smooth varieties by the rank of \( \alpha_{r-1} \) and apply Feix's
 result stratum by stratum, but to obtain a suitable hyperk\"ahler
thickening a more global approach is required.

The discussion in \( \S \)3 and \( \S \)4 suggests that we should consider first what the
analogue of the hypertoric variety \( Q_T \) might be when \( K \) is a symplectic or special
orthogonal group. As in \( \S \)2 let us first consider the case of \( K=\SO(n) \) when \( n=2r-1 \) is
odd. 

For the universal symplectic implosion in this case we considered symplectic quivers
\begin{equation}
  \label{eq:symplectic2}
  0  \stackrel{\alpha_0}{\rightarrow}
  \C \stackrel{\alpha_1}{\rightarrow}
  \C^2 \stackrel{\alpha_2}{\rightarrow} \dots
  \stackrel{\alpha_{r-2}}{\rightarrow} \C^{r-1}
  \stackrel{\alpha_{r-1}}{\rightarrow} 
   \C^n
\end{equation}  
and imposed the constraint 
\( \alpha_{r-1}^t J \alpha_{r-1} =0 \); we then took the symplectic quotient by
\( H_r = \prod_{j=1}^{r-1} \SU(j) \) with respect to the standard K\"ahler
structure on \( (\C^{j-1})^* \otimes \C^j \) for \( j \leqslant r-1 \) and the K\"ahler
structure induced by \( J \) on \( (\C^{r-1})^* \otimes \C^n \cong (\C^n)^{r-1} \).
We saw that there is a natural map to this symplectic quotient from the
toric variety given by the symplectic quotient of the space of symplectic quivers 
as above where each map \( \alpha_j \) has the form 
\begin{equation*}
  \alpha_j =
  \begin{pmatrix}
    0 & 0 & \dots &  0 \\
    \nu_1^j & 0 & \dots &  0 \\
    0 & \nu_2^j & \dots & 0 \\
    \vdots & \vdots & \ddots & \vdots \\
    0 & 0 &\dots &  \nu^j_j
  \end{pmatrix}
\end{equation*}
for \( j<r-1 \) and 
\begin{equation*}
\alpha_{r-1} = 
  \begin{pmatrix}
    0 & 0 & \dots & 0 \\
    \vdots & \vdots & & \vdots \\
    0 & 0 & \dots &  0 \\
    \nu_1^{r-1} & 0 & \dots &  0 \\
    0 & \nu_2^{r-1} & \dots & 0 \\
    \vdots & \vdots & \ddots & \vdots \\
    0 & 0 &\dots &  \nu^{r-1}_{r-1}
  \end{pmatrix}
\end{equation*}
with \( \nu_i^j \in \C \) as at \eqref{quiverT}; notice that a quiver of this
form always satisfies the constraint \( \alpha_{r-1}^t J \alpha_{r-1} =0 \).

By analogy with this and with the hypertoric variety \( Q_T \) described in
\( \S \)3 for the case when \( K=\SU(n) \), we expect the hypertoric variety \( Q_T \)
for \( K = \SO(n) \) when \( n=2r-1 \) to be closely related to the hyperk\"ahler quotient by the
standard maximal torus \( T_{H_r} \) of \( H_r \) of the flat space \( M_T^{\SO(n)} \) given 
by quiver diagrams
\begin{equation} \label{MSOn}
  0 \stackrel[\beta_0]{\alpha_0}{\rightleftarrows}
  \C \stackrel[\beta_1]{\alpha_1}{\rightleftarrows}
  \C^{2}\stackrel[\beta_2]{\alpha_2}{\rightleftarrows}\dots
  \stackrel[\beta_{r-2}]{\alpha_{r-2}}{\rightleftarrows} \C^{{r-1}}
  \stackrel[\beta_{r-1}]{\alpha_{r-1}}{\rightleftarrows}  \C^n
\end{equation}
where the maps \( \alpha_j \) and \( \beta_j \) have the form
\begin{equation*}   
  \alpha_j =
  \begin{pmatrix}
    0 & 0 & \dots &  0 \\
    \nu_1^j & 0 & \dots &  0 \\
    0 & \nu_2^j & \dots & 0 \\
    \vdots & \vdots & \ddots & \vdots \\
    0 & 0 &\dots &  \nu^j_j
  \end{pmatrix}
\quad
\text{and}
\quad
  \beta_j =
  \begin{pmatrix}
    0 & \mu_1^j & 0 & \dots &0  & 0 \\
    0 & 0 & \mu_2^j &  \dots &0 & 0 \\
    \vdots & \vdots & \vdots & \ddots & \vdots & \vdots \\
    0 & 0 & 0 & \dots &\mu_{j-1}^j & 0 \\
    0 & 0 & 0 &\dots &0 & \mu^j_j
  \end{pmatrix}
\end{equation*}
if \( j<r-1 \) and
\begin{equation*}   
\alpha_{r-1} = 
  \begin{pmatrix}
    0 & 0 & \dots & 0 \\
    \vdots & \vdots & & \vdots \\
    0 & 0 & \dots &  0 \\
    \nu_1^{r-1} & 0 & \dots &  0 \\
    0 & \nu_2^{r-1} & \dots & 0 \\
    \vdots & \vdots & \ddots & \vdots \\
    0 & 0 &\dots &  \nu^{r-1}_{r-1}
  \end{pmatrix}
\end{equation*}
and 
\begin{equation*}
  \beta_{r-1} =
  \begin{pmatrix}
    0&\dots&    0 & \mu_1^{r-1} & 0 & \dots &0  & 0 \\
    0&\dots &   0 & 0 & \mu_2^{r-1} &  \dots &0 & 0 \\
    \vdots&  &  \vdots & \vdots & \vdots & \ddots & \vdots & \vdots \\
    0&\dots&    0 & 0 & 0 & \dots &\mu_{r-2}^{r-1} & 0 \\
    0&\dots & 0 & 0 & 0 &\dots &0 & \mu^{r-1}_{r-1}
  \end{pmatrix}
\end{equation*}
for some \( \nu^j_i, \mu^j_i \in \C \setminus\{ 0 \} \). Notice that
\begin{equation} \label{eqnJ} 
\alpha_{r-1}^t J \alpha_{r-1} = 0 = \beta_{r-1} J \beta_{r-1}^t
\end{equation} 
for any \( \alpha_{r-1} \) and \( \beta_{r-1} \) of this form.

Let \( M^{\SO(n)} \) be the flat hyperk\"ahler space given by
arbitrary quiver diagrams of the form \eqref{MSOn}, where the
hyperk\"ahler structure is induced by the standard hyperk\"ahler
structure on \( (\C^{j-1})^* \otimes \C^j \oplus (\C^j)^* \otimes
\C^{j-1} \cong \HH^{(j-1)j} \) for \( j \leqslant r-1 \) and the
hyperk\"ahler structure induced by \( J \) on \( (\C^{r-1})^* \otimes
\C^n \oplus (\C^n)^* \otimes \C^{r-1} \cong (\HH^n)^{r-1} \).  As in
the symplectic case discussed in \( \S \)2, the restriction to \(
M_T^{\SO(n)} \) of the hyperk\"ahler moment map for the action of \(
H_r \) coincides with the hyperk\"ahler moment map for the action of
\( T_{H_r} \) on \( M_T^{\SO(n)} \). Thus
\begin{equation*}
Q_T^{\SO(n)} = M_T^{\SO(n)} \hkq T_{H_r}
\end{equation*}
maps naturally to \( M^{\SO(n)} \hkq H_r \). 

By analogy with the discussion in \( \S \)4 we can consider the subset
of \( M^{\SO(n)} \hkq H_r \) which is the closure of the \( K_\C =
\SO(n,\C) \)-sweep of the image of \( Q_T^{\SO(n)} \). By \eqref{eqnJ}
this is contained in the closed subset defined by the \( K_\C
\)-invariant constraints
\begin{equation*}
\alpha_{r-1}^t J \alpha_{r-1} = 0 = \beta_{r-1} J \beta_{r-1}^t,
\end{equation*}
and this closed subset has the dimension expected of the universal
hyperk\"ahler implosion.  Thus we expect the subset of the
hyperk\"ahler quotient \( M^{\SO(n)} \hkq H_r \) defined by these
constraints to 
be closely related to the universal hyperk\"ahler implosion for \( K =
\SO(n) \) when \( n=2r-1 \) is odd.  Similarly we expect that
modifications of this construction as described in \( \S \)2 for the
universal symplectic implosion will be closely related
to 
the universal hyperk\"ahler implosion for the special orthogonal
groups \( K = \SO(n) \) when \( n \) is even and for the symplectic
groups. The following example, however, provides a warning against
over-optimism here.

\begin{example}
  Recall that \( \SO(3) = \SU(2)/\{\pm 1 \} \) and that the universal
  symplectic implosion for \( \SO(3) \) is \( \C^2/\{ \pm 1\} \),
  where \( \C^2 \) is the universal symplectic implosion for \( \SU(2)
  \). Moreover the universal hyperk\"ahler implosion for \( \SU(2) \)
  is \( \HH^2 \), given by quiver diagrams of the form
  \begin{equation} \label{h2}
    \C^{}\stackrel[\beta]{\alpha}{\rightleftarrows}
    \C^{2},
  \end{equation}
  (recall that the group \( H_1 \) here is trivial, as in
  Example~\ref{SO3symp}, so no quotienting occurs).  We thus expect
  the universal hyperk\"ahler implosion for \( \SO(3) \) to be \(
  \HH^2 /\{\pm 1 \} \).

  We can associate to any quiver~\eqref{h2} the quiver
  \begin{equation*}
    \C^{} \cong \mathrm{Sym}^2(\C) \stackrel[\beta_1]{\alpha_1}{\rightleftarrows}
    \mathrm{Sym}^2(\C^2) \cong \C^3
  \end{equation*}
  where \( \alpha_1 = \mathrm{Sym}^2(\alpha) \) and \( \beta_1 =
  \mathrm{Sym}^2(\beta) \) are the maps between \( \mathrm{Sym}^2(\C)
  \) and \( \mathrm{Sym}^2(\C^2) \) induced by \( \alpha \) and \(
  \beta \). This construction gives us a surjection from \( \HH^2 \)
  to the subvariety of the space \( M^{\SO(3)}\hkq H_1 = M^{\SO(3)} \)
  of quivers
  \begin{equation*}
    \C^{} \stackrel[\beta_1]{\alpha_1}{\rightleftarrows}
    \C^{3}
  \end{equation*}
  satisfying \( \alpha_1^t J \alpha_1 = 0 = \beta_1 J \beta_1^t \),
  but it gives an identification of this subvariety with the quotient
  of \( \HH^2 \) by \( \Z_2 \times \Z_2 \), not by \( \Z_2 = \{ \pm 1
  \} \).
\end{example}

\end{document}